\documentclass[authoryear]{elsarticle}
\usepackage{hyperref}
\pdfoutput=1
\usepackage{listings}
\usepackage{color}
\journal{Econometric and Statistics: (Part B: Statistics)}
\usepackage{lineno}
\modulolinenumbers[5]
\usepackage{amsmath}
\usepackage{amsfonts}
\usepackage{amssymb}
\usepackage{amsthm}
\newtheorem{thm}{Theorem}
\theoremstyle{definition}
\newtheorem*{defn}{Definition}
\date{\today}
\title{Asymptotics for Markov chain mixture detection\tnoteref{t1}}
\hypersetup{
 pdfauthor={michaels},
 pdftitle={Asymptotics for Markov chain mixture detection\tnoteref{t1}},
 pdfkeywords={},
 pdfsubject={},
 pdfcreator={Emacs 26.3 (Org mode 9.3.1)}, 
 pdflang={English}}
\begin{document}

\begin{frontmatter}
\tnotetext[t1]{R source code for simulations in section 3.1 available on journal website}
\author[mf]{Matthew Fitzpatrick}
\address[mf]{Financial Markets, Westpac Banking Corporation}
\ead{Mattfitzpatrick@westpac.com.au}
\author[ms]{Michael Stewart\corref{ca}}
\cortext[ca]{Corresponding author}
\address[ms]{School of Mathematics and Statistics, University of Sydney}
\ead{michael.stewart@sydney.edu.au}

\begin{abstract}
Sufficient conditions are provided under which the log-likelihood ratio test statistic
fails to have a limiting chi-squared distribution under the null hypothesis when testing between one and two components under a general two-component mixture model, but rather tends to infinity in probability. These conditions are verified when the component densities describe continuous-time, discrete-state-space Markov chains and the results are illustrated via a parametric bootstrap simulation on an analysis of the migrations over time of a set of corporate bonds ratings. The precise limiting distribution is derived in a simple case with two states, one of which is absorbing which leads to a right-censored exponential scale mixture model. In that case, when centred by a function growing logarithmically in the sample size, the statistic has a  limiting distribution of Gumbel extreme-value type rather than chi-squared.
\end{abstract}

\begin{keyword}
Markov chain, mixture model, asymptotics
\MSC[2020] 60J28\sep 62M07
\end{keyword}
\end{frontmatter}


\section{Introduction.}
\label{sec:orge7e3a70}

Finite mixture models provide an appealling middle ground between
simple but narrow parametric models and broad but complex
nonparametric models. We posit a finite number of simple parametric
subpopulations or ``components'' and model that each observation is
drawn from one of the components via a random labelling mechanism, but
the labels are not observed. This idea is the basis of model-based
clustering, for which a central and still challenging aspect is the
question: ``how many clusters?''. The simplest version of this model
selection problem may be framed as a test of the null hypothesis of
one component versus the alternative of two components.  A very
important aspect of such a procedure is how one would calibrate the
test, given that usual regularity conditions which guarantee
asymptotic chi-squared distributions for standard test statistics are
violated.

The present work is motivated by \cite{frydman05} where a
two-component mixture of Markov chains was used as a potential model
for changes in rating level for financial entities over time.  This
work itself was proposed as a generalization of the ``mover-stayer''
model, itself a special case of a two-component mixture of Markov
chains where one of the components represents no movement between
states at all (see \cite{frydman05} for references to applications of
the mover-stayer model). The methodology is applied to the migration
of bond ratings of 848 corporate bond issuers between 7 different
rating levels (Aaa, Aa, A, Baa, Ba, B and C) and 2 absorbing states
(``default'' and ``rating withdrawal''). A log-likelihood ratio test (LRT)
between a 1-component and 2-component mixture is performed and
crucially an approximate p-value is computed assuming a limiting
chi-squared distribution under the null hypothesis.

\cite{hartigan85} was the first to point out the non-standard limiting
behaviour of the LRT statistic \(\Lambda_{n}\) for testing between one
and two components in a simple normal location mixture model where
there is no \emph{a priori} restriction on the possible parameter
values. He showed that \(\Lambda_n\) does not have a limiting
chi-squared distribution under the null hypothesis, but rather tends
to infinity in probability, conjecturing the rate to be \(\log\log n\)
where \(n\) is the sample size.  Using results in \cite{bickelchernoff93},
\cite{liushao04} verified Hartigan's conjecture by showing that
\(\Lambda_n-\log\log n\) has a limiting Gumbel (extreme-value)
distribution. \cite{liu03} showed the analogous result for a gamma scale
mixture model which includes an exponential scale mixture model as a
special case.

\cite{fukumizu03} generalized the result of \cite{hartigan85} to a very
general class of locally conic models as introduced by
\cite{dacunha97}. These models have a parametrization similar to
Euclidean polar coordinates. A fixed distribution corresponds to the
origin, while other members of the model are parameterised by two
parameters: a non-negative ``distance'' away from the origin and a
``direction''; the set of directions can be quite general. In particular
these models allow for a non-unique parametrisation of the ``origin'',
since setting the distance equal to zero, for any direction, returns
the null hypothesis. This nonidentifiability violates conventional
regularity conditions and leads to the non-standard asymptotics seen
in the simple mixture model tests cited above.

If one may reasonably restrict the parameter values of the mixture
components to a compact subset of the parameter space, different
asymptotic behaviour is obtained. Notable results for such cases were
obtained by, among others, \cite{ghoshsen85}, \cite{garel01,garel05},
\cite{dacunha97}, \cite{gassiat02}, \cite{liushao03}, \cite{azais06}, \cite{azais09}.  In these
cases \(\Lambda_n\) no longer tends to infinity in probability under the null
hypothesis, rather its limiting distribution is that of the square of
the maximum of a Gaussian process indexed by the (compact) parameter
space, or the difference of such maxima if the null hypothesis is
composite. While this restriction somewhat facilitates analysis,
justifying it in practice is not always possible. 
This issue is discussed further in section \ref{sec:orgabf1e47} below.

The non-regular limiting behaviour of maximised log-likelihoods when
testing between one and two mixture components leads to similar
non-standard asymptotics in cases  where up to three or more mixture
components are being considered. \cite{hui15} show that AIC (Akaike's Information Criterion) tends to overestimate the number of components, while Schwarz's BIC (Bayesian Information Criterion, and their novel adjusted version of AIC which they call ``AIC\textsubscript{mix}'' perform well. This somewhat mirrors the case in regular regression models where AIC tends to do better at prediction, while BIC tends to do better at model selection (see for example \cite{yang05}). However, the non-identifiability mentioned above introduces a subtle problem for these methods, which is further discussed in subsection \ref{sec:org2096391}.

The present paper contains two main results. The first provides a
convenient sufficient condition implying the result of \cite{fukumizu03}
in the particular setting of testing one versus two components in a
parametric mixture model. We derive this and apply it to the Markov
chain mixture model of \cite{frydman05}, showing that in that case the
LRT statistic tends to infinity in probability. Our second result
derives the precise limiting distribution in the simplest version of
this model where the Markov chain only has two states, one of which is
absorbing and all observations start in the nonabsorbing state. When
there is a finite observation window this leads to a censored version
of the exponential scale mixture model from \cite{liu03}; as in that
case we have \(\Lambda_n-\log\log n\) tending to a Gumbel distribution.

The structure of the paper is as follows: in section \ref{sec:org67d1a7f} we specify the general test of homogeneity in a
two-component mixture model and provide our first main result which
gives sufficient conditions for the LRT statistic to tend to infinity
in probability. In section \ref{sec:org9ae43fb}
we show how these conditons may be verified when the mixture
components are continuous-time Markov chains, which may have an
absorbing state. In section \ref{sec:org3066827} we indicate how the exact limiting
distribution of the LRT statistic may be derived in the simplest case
where the Markov chain has two states, one of which is absorbing and
all observations start in the non-absorbing state, so that times until
absorption are censored exponential random variables. Section
\ref{sec:org71ab6fd} has some concluding remarks and recommendations.

\section{Formulation of the general problem}
\label{sec:org67d1a7f}
\subsection{Two-component mixture models; tests of homogeneity}
\label{sec:orgf0a8485}
We are interested in the behaviour of the LRT when testing homogeneity
in a two-component mixture model and the null hypothesis of
homogeneity is true. We require that the component distributions are
taken from a parametric family which is suitably regular at the true
distribution. To be specific, suppose we have a parametric family
\(\mathcal{F}=\{f(\cdot;\theta,\gamma)\colon \theta\in \Theta,
\gamma\in\Gamma\}\) of densities with respect to a \(\sigma\)-finite
measure on a measurable space.
For our first main
result the family is otherwise quite general although in section
\ref{sec:org9ae43fb} we shall specialise to the
case where each \(f(\cdot;{\theta},\gamma)\) describes a continuous-time
Markov chain.

We shall specify regularity in terms of the limiting behaviour of
certain test statistics.  Suppose that our observations are modelled
as independent random elements
\(\boldsymbol{X}_1,\ldots,\boldsymbol{X}_n\) with common density some
\(f(\cdot;\theta_0,\gamma_0)\in \mathcal{F}\).  Define the one-component
log-likelihood as
\begin{align*}
L_n^{(1)}(\theta,\gamma) = \sum_{i=1}^n \log f(\boldsymbol{X}_i;{\theta},\gamma )\,.
\end{align*}

\begin{defn}[Regular Point]
We say that the density \(f(\cdot;\theta_0,\gamma_0)\) is a \emph{regular
point} of the family \(\mathcal{F}\) if the LRT statistic for testing
the simple hypothesis \((\theta,\gamma)=(\theta_0,\gamma_0)\) within
\(\mathcal{F}\), given by
\begin{align}
\label{eq:org3fbf7c8}
\Lambda_n^{(1)} =  \Lambda_n^{(1)}(\theta_0,\gamma_0)= \sup_{\substack{\theta\in\Theta\\ \gamma\in\Gamma}} L_n^{(1)}(\theta,\gamma) - L_n^{(1)}(\theta_0,\gamma_0)=O_p(1)\,, 
\end{align}
as \(n\to\infty\). 
\end{defn}
Here we use \(O_p(1)\) to denote a generic sequence of random variables
that remain bounded in probability.  Conditions implying (\ref{eq:org3fbf7c8})
are given in many places in the literature. Special cases include
where \(\Lambda_n^{(1)}\) is asymptotically \(\chi^2\)
\citep[see][Theorem~12.4.2]{lehmannromano05} or a mixture of
\(\chi^2\)s, \citep[see][]{chernoff54}.

Define the two-component mixture density
\begin{align*}
g(\cdot;p,\eta,\theta,\gamma) = (1-p) f(\cdot;\eta,\gamma) + p f(\cdot;\theta,\gamma)\,,
\end{align*}
where \(\eta,\theta\in\Theta\) are potentially different across
different components, \(\gamma\in\Gamma\) is the same across different
components and \(0\leq p\leq 1\). We consider the following two
different versions of a test between 1 and 2 components: the
simple-null-hypothesis version
\begin{align*}
H_0^s\colon &\boldsymbol{X}_1 \sim f(\cdot;{\eta_0},\gamma_0)\ \ \ \ \ \ \ \text{ for  known $\eta_0\in\Theta$, $\gamma_0\in\Gamma$ vs. }\\
H_1^{s}\colon &\boldsymbol{X}_1\sim g(\cdot;p,\eta_0,\theta,\gamma_0)\ \text{ for unknown $0<p<1$, $\theta\in\Theta$ and $\theta\neq \eta_0$;}
\intertext{and the composite-null-hypothesis version}
H_0^c\colon &\boldsymbol{X}_1 \sim f(\cdot;{\eta},\gamma)\ \ \ \  \ \ \ \text{ for $\eta\in\Theta$, $\gamma\in\Gamma$  vs.}\\
H_1^c\colon &\boldsymbol{X}_1\sim g(\cdot;{p,\eta,\theta,\gamma})\  \text{ for $0<p<1$,  $\theta\in\Theta$, $\theta\neq \eta$, }
\end{align*}
and all parameters are considered unknown.

Define the two-component log-likelihood as
and 
\begin{align*}
L_n^{(2)}(p,\eta,\theta,\gamma) = 
\sum_{i=1}^n \log g(\boldsymbol{X}_i;p,\eta,\theta,\gamma)\,.
\end{align*}
The LRT statistic for testing \(H_0^s\) versus \(H_1^s\) is
\begin{align}
\label{eq:orge649497}
\Lambda_n^s & = 
\sup_{\substack{0\leq p\leq 1\\ \theta\in\Theta}} 
L_n^{(2)}(p,\eta_0,\theta,\gamma_0) - L_n^{(1)}(\eta_0,\gamma_0)\,,
\end{align}
while the LRT statistic for testing \(H_0^c\) versus \(H_1^c\) 
is
\begin{align*}
\Lambda_n^c & = 
\sup_{\substack{0\leq p\leq 1\\ \eta,\theta\in\Theta\\ \gamma\in\Gamma}}
 L_n^{(2)}(p,\eta,\theta,\gamma) - \sup_{\eta\in\Theta}L_n^{(1)}(\eta)\,.
\end{align*}
We shall provide conditions below under which \(\Lambda_n^s\stackrel
P\to\infty\), that is for each \(0<M<\infty\),
\begin{align}
\label{eq:orga184d7a}
P \left\{ \Lambda_n^s\leq M \right\}\to 0\,.
\end{align}
We shall also be interested in the rate at which \(\Lambda_n^s\) tends
to infinity in probability.
\begin{defn}[Asymptotic Bounds; Rate of Divergence]
We say that a sequence \(r_n\to\infty\) provides an \emph{asymptotic upper
bound} to \(\Lambda_n^s\) if
\begin{align*}
\Lambda_n^s/r_{n} = O_p(1)\,,
\end{align*}
and an \emph{asymptotic lower bound} to \(\Lambda_n^s\) if
\begin{align*}
r_n/\Lambda_n^s = O_p(1)\,.
\end{align*}
We say that \(\Lambda_n^s\) \emph{diverges to infinity in probability at
rate} \(r_n\) if \(r_n\) is both an upper and lower asymptotic bound to
\(\Lambda_n^s\).
\end{defn}
Our regularity assumption means that if \(r_n\to\infty\) is an
asymptotic lower bound for \(\Lambda_n^s\) it is also an asymptotic
lower bound to \(\Lambda_n^c\). To see this, note that if
\(f(\cdot;\eta_0,\gamma_0)\) is a regular point of \(\mathcal{F}\),
\begin{align}
\label{eq:org5cf1355}
\Lambda_n^c \geq 
 \sup_{0\leq p\leq 1, \theta\in\Theta}
 L_n^{(2)}(p,\eta_0,\theta,\gamma_0)
 -  \sup_{\eta\in\Theta}L_n^{(1)}(\eta)
= \Lambda_n^s + O_p(1)\,.
\end{align}
We shall henceforth focus on \(\Lambda_n^s\), noting that results we
obtain are also applicable to \(\Lambda_n^c\), via (\ref{eq:org5cf1355}).

The profile likelihood
\begin{align*}
\Lambda_n^{p}(\theta)
=\sup_{0\leq p\leq 1}L_n^{(2)}(p,\eta_0,\theta,\gamma_0)\,,
\end{align*}
is the LRT statistic for testing \(H_0^s\) within the one-dimensional
submodel 
\begin{align}
\label{eq:orgc1198ff}
\{g(\cdot;p,\eta_0,\theta,\gamma_0)\colon 0\leq p\leq 1\}\,,
\end{align}
where the
parameter corresponding to the second component \(\theta\neq \eta_0\) is regarded as known.

Under certain conditions \(\Lambda_n^p(\theta)\) is, under \(H_0^s\),
asymptotically equivalent to half the square of the positive part of
the standardised scores with respect to \(p\) at \(p=0\). To make this precise, define the density ratio as
\begin{align*}
\nonumber r(\cdot;{\theta})&= f(\cdot;{\theta},\gamma_0) /f(\cdot;{\eta_0},\gamma_0)\,,
\intertext{and note that}
\left.\frac{\partial L_n^{(2)}(p,\eta_0,\theta,\gamma_0)}{\partial p}\right|_{p=0} &= \sum_{i=1}^n \left[ r(\boldsymbol{X}_i;{\theta})-1 \right]\,.
\end{align*}
Define  also
\begin{align*}
v_{\theta}&= E_0 \left[ r(\boldsymbol{X}_1;\theta)^2 \right]-1 = \mathrm{Var}_0 \left[ r(\boldsymbol{X}_1;{\theta}) \right]\,,
\intertext{and}
\Theta_v&=\left\{ \theta \colon v_{\theta} <\infty \right\}\,.
\end{align*}
Here  the zero subscript indicates expectation and variance when \(\boldsymbol{X}_1\sim f(\cdot;\theta_0,\gamma_0)\).
For \(\theta\in\Theta_v\), the standardised scores are
\begin{align}
\label{eq:orgb67bcac}
S_n(\theta)&= (n v_{\theta})^{-1/2}\sum_{i=1}^n \left[ r_{\theta} (\boldsymbol{X}_i)-1 \right]\,.
\end{align}

\begin{defn}[Asymptotically Normal Submodel]
We say the one-dimensional submodel (\ref{eq:orgc1198ff}) is asymptotically normal at \(p=0\) if 
 for each fixed \(\theta\in\Theta_v\) we have
\begin{align}
\label{eq:org3ed757f}
2 \Lambda_n^p(\theta) &= \max \left[ 0,S_n(\theta) \right]^2 + o_p(1)\,.
\end{align}
\end{defn}

The asymptotic equivalence (\ref{eq:org3ed757f}) relating to models where the
true value is on the boundary of the parameter space was first
discussed in \cite{chernoff54} and regularity conditions implying
(\ref{eq:org3ed757f}) are given there. Weaker conditions specific to locally
conic models were given in \cite{dacunha97} and \cite{fukumizu03}. 

\begin{thm}
Suppose that for each fixed \(\theta\in\Theta\cap\Theta_v\), (\ref{eq:org3ed757f}) holds and
that there exists a sequence \(\{\theta_i\}\) with each
\(\theta_i\in\Theta\cap \Theta_v\) such that
\begin{align}
\label{eq:org4115347}
v_{\theta_i} \to \infty\,.
\end{align}
Then \(\Lambda_n^s\stackrel P\to \infty\).
\label{org2cd1e3c}
\end{thm}

\begin{proof}
We first show that if \(\boldsymbol{X}\sim f(\cdot;\theta_0)\) then the sequence of random variables 
\begin{align}
\label{eq:org6a8f1b9}
v_{\theta_i}^{-1/2} \left[ r(\boldsymbol{X};{\theta_i})-1 \right]\stackrel P\to0\,.
\end{align}
Fix \(\varepsilon>0\). Then for large enough \(i\), \(-\varepsilon\leq -v_{\theta_i}^{-1/2}\leq v_{\theta_i}^{-1/2} \left[ r(\boldsymbol{X};{\theta_i})-1 \right]\) and so 
\begin{align*}
\left\{  \left| v_{\theta_i}^{-1/2} \left[ r(\boldsymbol{X};{\theta_i})-1 \right] \right| > \varepsilon\right\}
=\left\{   v_{\theta_i}^{-1/2} \left[ r(\boldsymbol{X};{\theta_i})-1 \right]  >\varepsilon\right\}
= \left\{r(\boldsymbol{X};{\theta_i}) > 1 + v_{\theta_i}^{1/2}\varepsilon  \right\}\,.
\end{align*}
However by Markov's inequality, since \(E_0 \left[ r(\boldsymbol{X};{\theta_i}) \right]\equiv1\), 
\begin{align*}
P_0 \left\{r(\boldsymbol{X};{\theta_i}) > 1 + v_{\theta_i}^{1/2}\varepsilon  \right\} \leq \frac{1}{1 + v_{\theta_i}^{1/2}\varepsilon}\to0\,,
\end{align*}
and thus
(\ref{eq:org6a8f1b9}) holds. The proof is completed by using the same argument as Theorem 1 of \cite{fukumizu03}.
\end{proof}

\subsection{Comparison with results obtained under a compactness assumption}
\label{sec:orgabf1e47}
For each \(\theta\in \Theta\cap \Theta_{v}\), \(S_n(\theta)\) given by (\ref{eq:orgb67bcac}) is
 asymptotically standard normal and so \(2\Lambda_n^p(\theta)\) has asymptotic
 distribution of \(\max \left[ 0,Z \right]^2\) where \(Z\) is standard
 normal. If the parameter space \(\Theta\) is not too large in a certain
 sense, the empirical process \(\{S_n(\theta) \colon \theta \in \Theta_v \cap \Theta\}\) converges in
 distribution to a tight Gaussian process indexed by \(\Theta\cap \Theta_v\) (that
 part of the parameter space \(\Theta\) for which the likelihood ratio
 \(r_{\theta}(\cdot)\) is square integrable). This is implied in those results mentioned in the Introduction making compactness assumptions, for example \cite{dacunha97,garel01,garel05} and \cite{liushao03}.

In such cases, the class of functions indexing the score process
\(\{S_n(\theta)\colon \theta \in \Theta_v \cap \Theta\}\) is a Donsker class and the rich literature
on empirical processes, established by many authors including
\cite{vanderwellner96} and \cite{geer00}, can be brought to bear
on the problem. In particular the maximum of the approximating
Gaussian process is bounded in probability, and under further
conditions so too is \(\Lambda_n^s\). However, in such cases the condition
(\ref{eq:org4115347}) does not hold.

When  (\ref{eq:org4115347}) does hold, the parameter space indexing the
empirical process is ``too large'' so that a tight Gaussian process with the
same covariance structure does not exist. In that case quite different methods are
needed to analyse the limiting behaviour of \(\Lambda_n^s\).
Examples where \(\Theta\) is one-dimensional  are provided by \cite{liu03}, \cite{liushao04} and our Theorem \ref{org71c58fc} below. In each of these, 
for a certain increasing sequence of compact subsets (intervals) of the parameter space \(\{\Theta_n\}\) with each \(\Theta_n\subset\Theta_{n+1}\) and \(\cup_{n=1}^{\infty}\Theta_n = \Theta\),  the score process \(\{S_n(\theta)\colon\theta\in\Theta_n\}\) may be approximated by a Gaussian process \(\{G_n(\theta)\colon \theta\in \Theta_n\}\) with the same covariance structure. Extremal theory for Gaussian processes (see for instance \cite{cramerleadbetter67}, \cite{leadbetter83} or \cite{piterbarg96}) then implies that the sequence of random variables \(\max_{\theta\in \Theta_n} \left[ 0,G_n(\theta) \right]^2\) tends to infinity in probability but when suitably normalised has a limiting Gumbel distribution. The approximation is accurate enough that with the same normalisation,
\(2\Lambda_n^s\)  has the same limiting distribution.

It is interesting to point out that 
\citet[Theorem 4]{azais06} showed that in general cases of this kind,
for \emph{any} sequence \(\{\Theta_n\}\) of increasing subsets, the corresponding Gaussian process functional \(\max_{\theta\in\Theta_n} \left[ 0,G_n(\theta) \right]^2\) has a limiting Gumbel distribution (once appropriately normalised), however they did not identify the particular rate which gives the limiting distribution of \(\Lambda_n^s\), rather it was used to show that under certain local alternatives the test has no limting power. 

\section{Mixtures of continuous-time Markov chains}
\label{sec:org9ae43fb}
In this section we focus on the case where each component density \(f(\cdot;\theta,\gamma)\) describes a stationary, continuous-time Markov chain observed over a finite time window \((0,T]\) with a finite state space \(\mathcal{W}=\{1,\ldots,w\}\). We may parametrise such a process via 3 parameters:
\begin{itemize}
\item \({\alpha}=\{\alpha(j)\colon j=1,\ldots,w\}\), the vector of probabilities describing the distribution of the initial state;
\item \({\beta}=\{\beta(j)=\beta_j\colon j=1,\ldots,w\}\), the vector of rates determining the average time the process stays in each state;
\item \({\gamma}=\{\gamma(j,k)\colon j,k=1,\ldots,w\}\), the matrix (with zeroes on the diagonal) of transition probabilities governing movements between states.
\end{itemize}
Following \cite{frydman05} we consider two-component mixtures of such processes, allowing for different rates \(\beta\) between the mixture components, but given that a transition occurs, the probabilities of each transition \(\gamma\) are the same across all mixture components. Thus in the notation of the previous section, the parameter which has possibly different values across the two mixture components is \(\theta=(\alpha,\beta)\)  and \(\gamma\) retains its meaning as the parameter which has the same value across both mixture components.

A typical sample point for a realisation of such a process is a step-function \(\{\boldsymbol{X}(t)\colon 0\leq t\leq T\}\) whose jump heights and jump locations may be characterised by a sequence of \(m+1\) pairs \(\{(z_0,t_0),(z_1,t_1),\ldots,(z_m,t_m)\}\) for \(m\geq0\):
\begin{itemize}
\item the process spends time \(t_0\leq T\) in the initial state \(z_0\);
\item if \(t_0<T\), it transitions to state \(z_1\neq z_0\) and spends time \(t_1\leq T-t_0\) there,
\end{itemize}
and so on until \(t_0+\ldots+t_m=T\).

 To derive the appropriate density function we may interpret the
\(z_j\)s as the first \(m+1\) observations on a discrete-time Markov
chain \(Z_0,Z_1,\ldots\) with transition matrix \(\gamma\) and given
\(Z_0=z_0,Z_1=z_1,\ldots\), conditionally independent random variables
\(\tilde{T}_0,\tilde{T}_1,\ldots\), with \(\tilde{T}_j\) exponential with
rate \(\beta_{z_j}\). The (random) number of transitions \(M\) is determined
via
\begin{align*}
M = \min \left\{ m \colon \sum_{j=0}^m \tilde{T}_j \geq T \right\}\,,
\end{align*}
and the actual sojourn times in each state are given by
\begin{align*}
T_j&=\tilde{T}_j\ \text{ for $j<M$}\\
\intertext{and}
T_M&= T-\sum_{j=0}^{M-1}T_{j}\,.
\end{align*}
Write \(m\) for the observed value of \(M\),
\begin{align}
\label{eq:orgb04b747}
\tau_j = \sum_{i=0}^{m} t_i 1 \left\{ z_i=j \right\}\,,
\end{align}
for the total time spent in state \(j\) and if \(m>0\), 
\begin{align}
\label{eq:org99da3ff}
n_{jk} = \sum_{i=1}^m 1 \left\{z_{i-1}=j,z_i=k  \right\}\,,
\end{align}
for the total number of transitions from state \(j\) to state \(k\).
The value of the density function at the point \(\left\{ \boldsymbol{x}(\cdot) \right\}\) 
with initial state \(z_0\) and sufficient statistics \(\{\tau_j\}\) and \(\{n_{jk}\}\) 
is then given by
\begin{align}
\label{eq:org9bdef26}
&f(\boldsymbol{x}(\cdot);(\alpha,\beta),\gamma)\\&
 = 
\begin{cases}
\alpha(z_0) e^{-\beta({z_0})T} & \text{ if $m=0$,}\\
\displaystyle{\alpha(z_0) \prod_{j=1}^{w} e^{-\beta_j \tau_j} \prod_{\substack{k=1\\\gamma_{j,k}>0}}^w \beta_j^{n_{jk}} \gamma_{j,k}^{n_{jk}}} & \text{ if $m>0$ and $\sum_{j=0}^m t_j=T$ and } \\
0 & \text{ otherwise.}
\end{cases}
\end{align}
This is precisely the form of the density as given in
\citet[Theorem~3.2]{albert62}, re-expressed in terms of the sufficient
statistics (\ref{eq:orgb04b747}) and (\ref{eq:org99da3ff}); see also that paper for a
description of the dominating measure associated with this
density. The restriction of the second product to indices \(k\) for
which \(\gamma_{j,k}>0\) accounts for both the zero diagonal elements of
\(\gamma\) as well as the case when any of the states are absorbing.

We wish to apply Theorem \ref{org2cd1e3c} from the previous section to the version
of \(\Lambda_n^s\) obtained when the class of densities \(\mathcal{F}\) is
given by (\ref{eq:org9bdef26}). We write the two-component density as
\begin{align*}
g(\,\cdot\,;p,(\alpha_0,\beta_{0}),(\alpha_1,\beta_1),\gamma_0)=(1-p) f(\cdot;(\alpha_0,\beta_0),\gamma_0)
+p f(\cdot;(\alpha_1,\beta_1),\gamma_0)\,.
\end{align*}
Note that the transition probability matrix \(\gamma_{0}\) is the same for both and is considered known.

We first show that condition (\ref{eq:org4115347}) from Theorem \ref{org2cd1e3c} is satisfied. 
Consider a sequence of parameter values \(\{(\alpha_i,\beta_i)\}\)
defined by \(\alpha_i\equiv \alpha_0\) and \(\beta_i = c_i \beta_0\) for a
positive real sequence \(c_i\to\infty\).  It is
straightfoward to see that the distributions corresponding to
\(((\alpha_i,\beta_i),\gamma_0)\) and \(((\alpha_0,\beta_0),\gamma_0)\)
have the same support. Consequently, the corresponding score variance
\(v_{\theta_i}=v_{(\alpha_i,\beta_i)}\) is the integral with respect to
the dominating measure of the function
\begin{align*}
\frac{\left[ f(\cdot; (\alpha_0,c_i\beta_0),\gamma_0)
	  \right]^2}{f(\cdot;(\alpha_0,\beta_0),\gamma_0)}-1\,.
\end{align*}
The first term evaluated at a sample point \(\boldsymbol{x}(\cdot)\) with intial state \(z_0\), \(m\) transitions 
and sufficient statistics \(\{\tau_j\}\) and \(\left\{ n_{jk} \right\}\) is equal to
\begin{align}
\label{eq:orge149251}
\nonumber &\frac{\left[ f(\boldsymbol{x}(\cdot); (\alpha_0,c_i\beta_0),\gamma_0)
	  \right]^2}{f(\boldsymbol{x}(\cdot);(\alpha_0,\beta_0),\gamma_0)}\\
\nonumber &= 
\begin{cases}
\alpha_0(z_0)e^{-(2c_i-1)\beta_0(z_0)} & \text{ for $m=0$,}\\
\displaystyle{\alpha_0(z_0)\prod_{j=1}^w e^{-(2c_i-1)\beta_0(j) \tau_j} \prod_{\substack{k=1\\ \gamma_0(j,k)>0}}^w  \left[ c_i^2 \beta_0(j)\right]^{n_{jk}} \gamma_0(j,k)^{n_{jk}}} & \text{ otherwise,}
\end{cases}\\
 &= \left( \frac{c_i^2}{2c_i-1}\right)^m f(\boldsymbol{x}(\cdot);\alpha_0,(2c_i-1)\beta_0,\gamma_0)\\
\intertext{and}
\nonumber &\geq 
 \left( \frac{c_i^2}{2c_i-1}\right) 1 \left\{m>0\right\} f(\boldsymbol{x}(\cdot);\alpha_0,(2c_i-1)\beta_0,\gamma_0)\,.\\
\intertext{The integral of this last right-hand side is}
\nonumber & 
 \left( \frac{c_i^2}{2c_i-1}\right) \left[1- P_{((\alpha_0,(2c_i-1)\beta_0),\gamma_0)}(M=0)\right]\\
\nonumber &=
 \left( \frac{c_i^2}{2c_i-1}\right) \left[1- \sum_{j=1}^w \alpha_0(j)e^{-(2c_i-1)\beta_0(j)T}\right]\,.
\end{align}
As \(c_i\to\infty\) the first factor tends to infinity while the second tends to 1. 
Also note that for fixed \(c_i\) the integral of (\ref{eq:orge149251}) is always finite since it is 
\begin{align*}
E_{(\alpha_0,\beta_0),\gamma_0} \left\{ e^{M \left[ 2 \log c_i - \log(2c_i-1) \right]} \right\}\,,
\end{align*}
but \(E_{(\alpha_0,\beta_0),\gamma_0} \left( e^{tM} \right)<\infty\) for
all real \(t\) as \(M\) is stochastically smaller than what would be
obtained if all rates were equal to \((2c_i-1)\left[ \max_j \beta_0(j)
\right]\), in which case \(M\) would be a Poisson random variable. This
verifies (\ref{eq:org4115347}).

We do not verify condition (\ref{eq:org3ed757f}) but instead refer the reader to
\citet[Section 4.3.3]{fitzpatrick16} where it is demonstrated that the
asymptotic normality conditions of \cite{fukumizu03} are satisfied in the
present case.
\subsection{Frydman's log-likelihood ratio test}
\label{sec:orgf2debc7}
\cite{frydman05} fits a two-component mixture of continuous-time Markov
chains to some bonds ratings data and tests the null hypothesis that
the true distribution is a single Markov chain. The state space
consists of 9 states, one of which (``default'') is an absorbing state.
We note here that a different parametrisation of the two-component
Markov chain is used, however the model is indeed equivalent to the
one we have formulated here; see \citet[Section~4.2]{fitzpatrick16}
for explanatory details.  The log-likelihood ratio statistic took the
value 276.96, and was judged to be highly significant when compared to
the \(\chi^2_8\) distribution. However, our development above shows that
the \(\chi^2\) distribution is not appropriate; under the null
hypothesis of a single Markov chain our theory tells us that the LRT
statistic tends to infinity in probability. However, our results do
not give any indication of the rate of divergence (in the next section
we see that an asymptotic lower bound is \(\log\log n\), where \(n\) is
the sample size).

\begin{figure}[htbp]
\centering
\includegraphics[width=.9\linewidth]{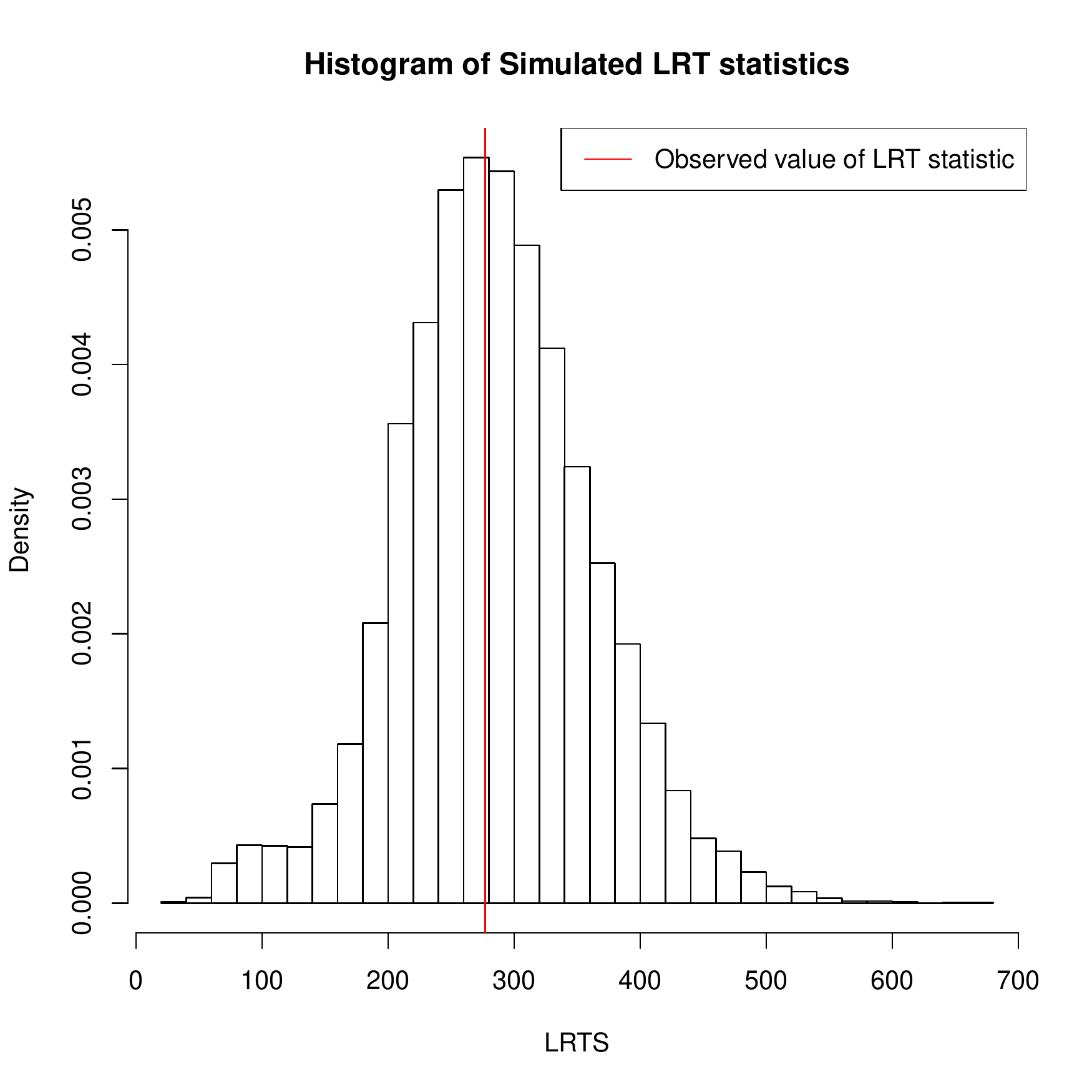}
\caption{\label{fig:org56ee73b}Histogram of LRT statistics obtained by simulating from the fitted one-component Markov chain from Frydman (2005). The approximate parametric bootstrap p-value is 0.528, suggesting the data is consistent with a single Markov chain as opposed to a mixture of two.}
\end{figure}

To better judge the significance of the observed value of the LRT
statistic, we performed a parametric bootstrap test of the same
hypothesis, obtaining an approximate p-value via simulation. The
parameters of the one-component fit were given in \cite{frydman05} so we
simulated from the one-component fit 10,000 times, fitting a one- and
two-component Markov chain mixture to each such sample and obtained
the LRT statistic each time. The proportion of simulated statistics
exceeding the observed value of 276.96 was 0.528 (see Figure \ref{fig:org56ee73b}). This indicates that
using the \(\chi^2_8\) as an approximate sampling distribution of the
LRT statistic leads to extreme false significance. The R source code for reproducing the parametric bootstrap simulation is supplied as supplementary material.

\subsection{Selecting the number of mixture components using standard information criteria}
\label{sec:org2096391}
Well-known model selection methods, particularly AIC, BIC and many other variants, all have the same general form: choose the model with the largest value of the difference
\begin{align*}
\Lambda_k - \pi_k\,,
\end{align*}
where \(\Lambda_k\) is the maximised log-likelihood over model \(k\) and \(\pi_k\) is a non-negative penalty which measures, in some sense, the ``complexity'' of model \(k\).

Any such method is determined by the the \emph{pairwise differences}
\(\pi_{\ell}-\pi_k\) for each \(k\neq \ell\). For many methods, including AIC and BIC,
this difference is a function of \(d\), the ``difference in number of
free parameters'' between the two models, alternatively the ``number of
parametric constraints'' imposed on the larger model to give the
smaller model. For AIC it is \(2d\), for BIC it is \(d \log n\), where \(n\)
is the sample size.  

These methods suffer from an operational problem
when trying to select the number of components in finite mixture
models, because smaller models with fewer components are not uniquely
identified within the larger models with more components.  For
example, consider a \(k\)-dimensional parametric family of densities
\(\{f(\cdot;\boldsymbol{\theta})\colon \boldsymbol{\theta}\in \mathbb{R}^k\}\) and define the
two-component mixture model
\begin{align*}
g(\cdot;p,\boldsymbol{\theta}) = (1-p) f(\cdot;\boldsymbol{\theta}_0) + p f(\cdot;\boldsymbol{\theta})\,,
\end{align*}
for some known reference value \(\boldsymbol{\theta}_0\). Take as the smaller
model the singleton \(\{f(\cdot;\boldsymbol{\theta}_0)\) and take as the larger
model \(\{g(\cdot;p,\boldsymbol{\theta})\colon 0\leq p\leq 1,\,\boldsymbol{\theta}\in
\mathbb{R}^k\}\).  The smaller model may be obtained from the larger
model by imposing the \(k\) constraints
\(\boldsymbol{\theta}=\boldsymbol{\theta}_0\), but it may also be obtained by
imposing the \emph{single} constraint \(p=0\). Thus the value of \(d\) here is
ambiguous, so it is not clear how penalties for AIC and BIC are
determined in this or similar cases.

We therefore caution against using information criteria of this form,
which only depend on this value \(d\), like AIC and BIC. More
complicated methods are needed which take into account the
non-standard behaviour of maximised likelihoods in these mixture
models, for example AIC\textsubscript{mix} introduced in \cite{hui15}. A more
detailed examination of this and similar methods is beyond the scope
of the current paper, where the focus is on choosing between one and
two components under the null hypothesis of one component. Once the
limiting null distribution of the log-likelihood ratio statistic is
determined, the limiting probability of selecting either model, using
any form of complexity penalty \emph{when the null hyothesis is true} can
also be determined.

\section{Exact asymptotics for detection of a  two-state Markov chain mixture}
\label{sec:org3066827}
We showed in the previous section that when fitting a mixture of two
Markov chains (even if one of the states is absorbing) the LRT
statistic for testing between one and two components tends to infinity
in probability. Our simulation also showed that the rate of divergence
is such that an observed statistic which seems highly significant when
compared to a \(\chi^2\) distribution may not be significant at all.

In this section we consider the simplest possible special case: where
there are only two states, one of which is absorbing and all observations
start in the non-absorbing state. Suppose that defaults occur with rate \(\theta\). 
Then the time to default has an
exponential distribution but if the observation window is limited to
an interval \((0,T]\) there is a positive probability of defaulting
after time \(T\).  Thus we may take each observation \(X_i=\min(Y_i,T)\)
where \(Y_i\) has exponential distribution with rate \(\theta\) so that
\(X_i\) has a right-censored exponential distribution satisfying
\begin{align*}
P_{\theta} (X_i\leq x) = 
\begin{cases}
0 & \text{ for $x\leq 0$,}\\
1  - e^{-\theta x} & \text{ for $0< x<T$ and}\\
1 & \text{ for $x=T$.} 
\end{cases}
\end{align*}
We may take as dominating measure the sum of Lebesgue measure on the
interval \((0,T)\) and counting measure on \(\{T\}\). The corresponding
density is then given by
\begin{align*}
f(x;\theta) = \theta e^{-\theta x} 1 \left\{ 0<x<T \right\} + e^{-\theta T} 1 \left\{ x=T \right\}\,.
\end{align*}
Without loss of generality we assume that \(X_1,\ldots,X_n\) are independent with common density
\(f(\cdot;1)\). The statistic (\ref{eq:orge649497}) becomes
\begin{align}
\label{eq:org547f18b}
\Lambda_n^s = \sup_{\substack{0\leq p\leq 1\\ \theta>0}} \sum_{i=1}^n \log \left\{ 1 + p \left[ \frac{f(X_i;\theta)}{f(X_i;1)}-1 \right] \right\}\,.
\end{align}

\begin{thm}
The statistic \(\Lambda_n^s\) given by (\ref{eq:org547f18b}) satisfies
\begin{align*}
\Lambda_n^s = \log\log n + G_n - \log 4\pi\,,
\end{align*}
where \(G_n\) has an asymptotic Gumbel distribution, that is
\begin{align}
\label{eq:orgbcd4bc6}
P\{G_n\leq x\} \to e^{-e^{-x}}\,.
\end{align}
\label{org71c58fc}
\end{thm}
The method of proof closely follows that of the corresponding result
in the one-sided, uncensored case, \citet[Theorem 2]{liu03} whose
steps we outline in subsection \ref{sec:org6479766}.  In subsection \ref{sec:org931429e},
we describe the modifications required for proving the corresponding
result in the two-sided censored case.

\subsection{The one-sided, uncensored case of \cite{liu03}.}
\label{sec:org6479766}
 \cite{liu03} derived the limiting null distribution of the statistic
\(\Lambda_n^s\) in the case of gamma scale mixtures (which includes
uncensored exponential scale mixtures), but the test was one-sided in
that the mean in the unknown second component was restricted to be
\emph{smaller} than the mean in the known, null component.  We follow the
same general approach for our censored case but with several notable
differences, in particular that we are able to relax the one-sided
restriction, mainly due to the censoring, as we shall see.

The main idea behind their method is to obtain a uniform version of
the approximation (\ref{eq:org3ed757f}) with careful control of the rate of the
error term. The main tool used to achieve this is the approximation (on a
suitable probability space) a version of the process \(\{S_n(\theta)\}\)
with a mean-zero Gaussian process \(\{H_n(\theta)\}\) which has the
property that the scale-transformed process \(\left\{H_n\left(e^s+
\frac{1}{2}\right)\right\}\) is  \emph{stationary},  satisfying
\begin{align*}
\rho(s,t) 
=E \left[ H_n\left(e^s+ \frac{1}{2}\right)H_n\left( e^t+\frac{1}{2}\right) \right]
=2\left(e^{\frac{s-t}{2}} + e^{\frac{t-s}{2}}\right)^{-1}\,,
\end{align*}
only depending on \((s,t)\) via \(|s-t|\); moreover as \(\Delta\to0\), 
\begin{align}
\label{eq:org8e6edc1}
\rho(t,t+\Delta)= 1 - \frac{\Delta^2}{8} + o \left( \Delta^2 \right)\,,
\end{align}
while as \(\Delta\to\infty\), 
\begin{align}
\label{eq:orgc286781}
\rho(t,t+\Delta)\sim 2 e^{-\frac{\Delta}2}\,.  
\end{align}
The short-range condition (\ref{eq:org8e6edc1}) and long-range condition
(\ref{eq:orgc286781}) mean that according to
\citet[Theorem~8.2.7]{leadbetter83}, for any \(L_n\to\infty\),
\begin{align}
\label{eq:org9a263e6}
\sup_{0\leq s\leq L_n} H_n \left( e^s+\frac{1}{2} \right) = \sqrt{2 \log L_n} + \frac{G_n - \log(4\pi)}{\sqrt{2 \log L_n}}\,,
\end{align}
where \(G_n\) denotes a generic random variable with an asymptotic Gumbel distribution, satisfying (\ref{eq:orgbcd4bc6}) above.
By stationarity, the same applies to any interval of length \(L_n\).

The remaining steps in their proof may be summarised as follows. Writing 
\begin{align*}
I_n = \left(\log n, n(\log n)^{-4}\right)\,,
\end{align*}
they show that
\begin{align*}
\sup_{\theta\in [1,\infty)\setminus  I_{n}} 
S_n(\theta) 
  &= O_p \left\{ (\log\log\log n)^{1/2} \right\}\,\\
\intertext{and}
\sup_{\theta\in I_n} 
\left| S_n(\theta)-H_n(\theta) \right| &= O_p \left\{ (\log n)^{-1} \right\}\,.
\end{align*}
These two results together with (\ref{eq:org9a263e6}) imply that with probability tending to 1, 
\begin{align*}
M_n&= \sup_{\theta\geq 1}S_n(\theta) = \sup_{\theta\in I_n} 
S_n(\theta) = \sup_{\theta\in I_n} 
H_n(\theta) + O_p \left\{ (\log n)^{-1} \right\}\,.
\end{align*}
This immediately yields
\begin{align*}
M_n = \sqrt{2 \log \log n} + \frac{G_n - \log 4\pi}{\sqrt{2 \log\log n}}\,,
\end{align*}
where again \(G_n\) satisfies (\ref{eq:orgbcd4bc6}), and the same is true replacing \(M_n\) with \(\sup_{\theta\in I_n}S_n(\theta)\).

The final steps involve showing that the profile likelihood \(\Lambda_n^p(\theta)\) satisfies
\begin{align*}
\sup_{\theta\in[1,\infty)\setminus I_n} 
\Lambda_n^p(\theta) &= O_p \left\{ (\log\log\log n)^{1/2} \right\}\\
\intertext{and}
\sup_{\theta\in I_n} 
\Lambda_n^p(\theta) &= \frac{1}{2}\left\{\sup_{\theta\in I_n} 
S_n(\theta)\right\}^{2} + o_p(1)\,.
\end{align*} 
These imply that
\begin{align*}
\Lambda_n^s = \sup_{\theta\geq 1}\Lambda_n^p(\theta)= \log\log n  + G_n - \log 4\pi\,,
\end{align*}
with \(G_n\) satisfying (\ref{eq:orgbcd4bc6}).

\subsection{The two-sided censored case}
\label{sec:org931429e}
\cite{liu03} only considered \(\theta\geq 1\) in their second mixture
component because the same strategy cannot be applied for
\(0<\theta<1\); in particular for \(0<\theta<\frac{1}{2}\) the score
process has infinite variance and so the Gaussian process
approximation fails. In short, the method fails if the second mixture
component has a mean which is too large. However in our case, the
censoring has the result of attenuating that problem and so it becomes
possible to approximate the score process \(\{S_n(\theta)\}\) with a
Gaussian process \(\{H_n(\theta)\}\) throughout the whole range
\(\theta>0\). The main modification required is to manage the fact that
this 
Gaussian process 
is such that using the same change
of scale,
\begin{align*}
\rho(t,t+\Delta) 
&=E \left[ H_n\left(e^t+ \frac{1}{2}\right)H_n\left( e^{t+\Delta}+\frac{1}{2}\right) \right]\\
\intertext{and}
&=\frac 2{e^{\frac{\Delta}{2}} + e^{-\frac{\Delta}{2}}} \, \frac{1 - e^{-T(e^t+e^{t+\Delta})}}{\sqrt{(1-e^{-2Te^t})(1-e^{-2Te^{t+\Delta}})}}
\end{align*}
which is not free of \(t\) as in the uncensored case. Thus the process
\(\left\{ H_n \left( e^t + \frac{1}{2} \right) \right\}\) is not
stationary. However, we show that it is \emph{locally stationary} in the
sense of \cite{berman74} and \cite{huesler90}, meaning that as
\(\Delta\to0\),
\begin{align*}
\rho(t,t+\Delta) = 1 - V(t)\Delta^2 + o \left( \Delta^2 \right)\,,
\end{align*}
for a  continuous positive function \(V(t)\) which is bounded away from zero and
infinity and the \(o \left( \Delta^2 \right)\) term is uniform in
\(t\) \citep[see][subsection 5.3.1 for details]{fitzpatrick16}. Under an analogue of the long-range condition (\ref{eq:orgc286781}),
namely
\begin{align*}
\sup_t \rho(t,t+\Delta) = o \left( \frac{1}{ \log \Delta}  \right)\,,
\end{align*}
as \(\Delta\to\infty\) (which is comfortably satisfied), a parallel
theory of extremes for locally stationary Gaussian processes
\citep[see][]{huesler90,huesler95} provides analogous limiting
results.

The rest of the proof follows the same steps as for
\citet[Theorem~2]{liu03}, with various technical modifications
needed. We do not provide all the details here, but refer the reader
to \citet[Chapter~5]{fitzpatrick16}.
\section{Concluding remarks}
\label{sec:org71ab6fd}
We have provided some general conditions which indicate when the
standard \(\chi^2\) approximation to the LRT statistic fails and may
lead to false significance; our simulations in subsection \ref{sec:orgf2debc7} show that this discrepancy can be
substantial. Our general advice is to always verify asymptotic
approximations to p-values by simulation whenever possible. The
precise limiting distribution we provide in section \ref{sec:org3066827} for the case of two
states indicates a very slow rate of divergence of \(\log\log n\),
however our simulation in section \ref{sec:orgf2debc7}
shows that for a larger number of states the rate of divergence may be
substantially faster. The rate is related to the complexity of the
family of standardised score functions and both the dimensionality and
effective range of maximisation (respectively 1 and \(L_n\approx\log n\)
in our example) of a corresponding approximating locally stationary
Gaussian process. The dimensionality reflects the number of states
while the effective range of maximisation reflects how the size of the
convex hull of the sample is transformed under the appropriate
local-stationarity-inducing transformation. Some indications of this
are given in some elementary univariate exponential family mixture
model examples in \cite{stewartrobinson03}; see also
\cite{ingster97,ingster01,ingster02}, \cite{donohojin04},
\cite{cai11,caiwu14,hallstewart05} and \cite{porterstewart20} for
more on rates of convergence of the LRT statistic under both the null
hypothesis and local alternatives in various two-component mixture models.

\bibliographystyle{model5-names}
\bibliography{Markov-asymptotics-Fitzpatrick-Stewart}

\begin{thebibliography}{37}
\expandafter\ifx\csname natexlab\endcsname\relax\def\natexlab#1{#1}\fi
\providecommand{\url}[1]{\texttt{#1}}
\providecommand{\href}[2]{#2}
\providecommand{\path}[1]{#1}
\providecommand{\DOIprefix}{doi:}
\providecommand{\ArXivprefix}{arXiv:}
\providecommand{\URLprefix}{URL: }
\providecommand{\Pubmedprefix}{pmid:}
\providecommand{\doi}[1]{\href{http://dx.doi.org/#1}{\path{#1}}}
\providecommand{\Pubmed}[1]{\href{pmid:#1}{\path{#1}}}
\providecommand{\bibinfo}[2]{#2}
\ifx\xfnm\relax \def\xfnm[#1]{\unskip,\space#1}\fi
\bibitem[{Albert(1962)}]{albert62}
\bibinfo{author}{Albert, A.} (\bibinfo{year}{1962}).
\newblock \bibinfo{title}{Estimating the infinitesimal generator of a
  continuous time, finite state {M}arkov process}.
\newblock {\it \bibinfo{journal}{Ann. Math. Statist.}\/},  {\it
  \bibinfo{volume}{33}\/}, \bibinfo{pages}{727--753}.
\bibitem[{Aza\"{\i}s et~al.(2006)Aza\"{\i}s, Gassiat \& Mercadier}]{azais06}
\bibinfo{author}{Aza\"{\i}s, J.-M.}, \bibinfo{author}{Gassiat, E.}, \&
  \bibinfo{author}{Mercadier, C.} (\bibinfo{year}{2006}).
\newblock \bibinfo{title}{Asymptotic distribution and local power of the
  log-likelihood ratio test for mixtures: bounded and unbounded cases}.
\newblock {\it \bibinfo{journal}{Bernoulli}\/},  {\it \bibinfo{volume}{12}\/},
  \bibinfo{pages}{775--799}. \URLprefix
  \url{https://doi.org/10.3150/bj/1161614946}.
  \DOIprefix\doi{10.3150/bj/1161614946}.
\bibitem[{Aza\"{\i}s et~al.(2009)Aza\"{\i}s, Gassiat \& Mercadier}]{azais09}
\bibinfo{author}{Aza\"{\i}s, J.-M.}, \bibinfo{author}{Gassiat, E.}, \&
  \bibinfo{author}{Mercadier, C.} (\bibinfo{year}{2009}).
\newblock \bibinfo{title}{The likelihood ratio test for general mixture models
  with or without structural parameter}.
\newblock {\it \bibinfo{journal}{ESAIM Probab. Stat.}\/},  {\it
  \bibinfo{volume}{13}\/}, \bibinfo{pages}{301--327}. \URLprefix
  \url{https://doi.org/10.1051/ps:2008010}. \DOIprefix\doi{10.1051/ps:2008010}.
\bibitem[{Berman(1974)}]{berman74}
\bibinfo{author}{Berman, S.~M.} (\bibinfo{year}{1974}).
\newblock \bibinfo{title}{Sojourns and extremes of {G}aussian processes}.
\newblock {\it \bibinfo{journal}{Ann. Probability}\/},  {\it
  \bibinfo{volume}{2}\/}, \bibinfo{pages}{999--1026}.
\bibitem[{Bickel \& Chernoff(1993)}]{bickelchernoff93}
\bibinfo{author}{Bickel, P.}, \& \bibinfo{author}{Chernoff, H.}
  (\bibinfo{year}{1993}).
\newblock \bibinfo{title}{Asymptotic distribution of the likelihood ratio
  statistic in a prototypical non regular problem}.
\newblock In \bibinfo{editor}{J.~K. Ghosh}, \bibinfo{editor}{S.~K. Mitra},
  \bibinfo{editor}{K.~R. Parthasararthy}, \& \bibinfo{editor}{B.~L.~S.
  Prakasa~Rao} (Eds.), {\it \bibinfo{booktitle}{Statistics and Probability: A
  Raghu Raj Bahadur Festschrift}\/} (pp. \bibinfo{pages}{83--96}).
\newblock \bibinfo{publisher}{Wiley Eastern Limited}.
\bibitem[{Cai et~al.(2011)Cai, Jeng \& Jin}]{cai11}
\bibinfo{author}{Cai, T.~T.}, \bibinfo{author}{Jeng, X.~J.}, \&
  \bibinfo{author}{Jin, J.} (\bibinfo{year}{2011}).
\newblock \bibinfo{title}{Optimal detection of heterogeneous and
  heteroscedastic mixtures}.
\newblock {\it \bibinfo{journal}{J. R. Stat. Soc. Ser. B Stat. Methodol.}\/},
  {\it \bibinfo{volume}{73}\/}, \bibinfo{pages}{629--662}.
\bibitem[{Cai \& Wu(2014)}]{caiwu14}
\bibinfo{author}{Cai, T.~T.}, \& \bibinfo{author}{Wu, Y.}
  (\bibinfo{year}{2014}).
\newblock \bibinfo{title}{Optimal detection of sparse mixtures against a given
  null distribution}.
\newblock {\it \bibinfo{journal}{IEEE Trans. Inform. Theory}\/},  {\it
  \bibinfo{volume}{60}\/}, \bibinfo{pages}{2217--2232}. \URLprefix
  \url{https://doi.org/10.1109/TIT.2014.2304295}.
  \DOIprefix\doi{10.1109/TIT.2014.2304295}.
\bibitem[{Chernoff(1954)}]{chernoff54}
\bibinfo{author}{Chernoff, H.} (\bibinfo{year}{1954}).
\newblock \bibinfo{title}{On the distribution of the likelihood ratio}.
\newblock {\it \bibinfo{journal}{Annals of Mathematical Statistics}\/},  {\it
  \bibinfo{volume}{25}\/}, \bibinfo{pages}{573--578}.
\bibitem[{Cram{\'e}r \& Leadbetter(1967)}]{cramerleadbetter67}
\bibinfo{author}{Cram{\'e}r, H.}, \& \bibinfo{author}{Leadbetter, M.~R.}
  (\bibinfo{year}{1967}).
\newblock {\it \bibinfo{title}{Stationary and related stochastic processes.
  {S}ample function properties and their applications}\/}.
\newblock \bibinfo{address}{New York}: \bibinfo{publisher}{John Wiley \& Sons
  Inc.}
\bibitem[{Dacunha-Castelle \& Gassiat(1997)}]{dacunha97}
\bibinfo{author}{Dacunha-Castelle, D.}, \& \bibinfo{author}{Gassiat, E.}
  (\bibinfo{year}{1997}).
\newblock \bibinfo{title}{Testing in locally conic models, and application to
  mixture models}.
\newblock {\it \bibinfo{journal}{ESAIM: Probability and Statistics}\/},  {\it
  \bibinfo{volume}{1}\/}, \bibinfo{pages}{285--317}.
\bibitem[{Donoho \& Jin(2004)}]{donohojin04}
\bibinfo{author}{Donoho, D.}, \& \bibinfo{author}{Jin, J.}
  (\bibinfo{year}{2004}).
\newblock \bibinfo{title}{Higher criticism for detecting sparse heterogeneous
  mixtures}.
\newblock {\it \bibinfo{journal}{Ann. Statist.}\/},  {\it
  \bibinfo{volume}{32}\/}, \bibinfo{pages}{962--994}.
\bibitem[{Fitzpatrick(2016)}]{fitzpatrick16}
\bibinfo{author}{Fitzpatrick, M.~A.} (\bibinfo{year}{2016}).
\newblock {\it \bibinfo{title}{Multi-regime models involving Markov chains}\/}.
\newblock Ph.D. thesis University of Sydney.
\newblock \URLprefix \url{http://hdl.handle.net/2123/14530}.
\bibitem[{Frydman(2005)}]{frydman05}
\bibinfo{author}{Frydman, H.} (\bibinfo{year}{2005}).
\newblock \bibinfo{title}{Estimation in the mixture of {M}arkov chains moving
  with different speeds}.
\newblock {\it \bibinfo{journal}{J. Amer. Statist. Assoc.}\/},  {\it
  \bibinfo{volume}{100}\/}, \bibinfo{pages}{1046--1053}. \URLprefix
  \url{http://dx.doi.org/10.1198/016214505000000024}.
  \DOIprefix\doi{10.1198/016214505000000024}.
\bibitem[{Fukumizu(2003)}]{fukumizu03}
\bibinfo{author}{Fukumizu, K.} (\bibinfo{year}{2003}).
\newblock \bibinfo{title}{Likelihood ratio of unidentifiable models and
  multilayer neural networks}.
\newblock {\it \bibinfo{journal}{Ann. Statist.}\/},  {\it
  \bibinfo{volume}{31}\/}, \bibinfo{pages}{833--851}. \URLprefix
  \url{https://doi.org/10.1214/aos/1056562464}.
  \DOIprefix\doi{10.1214/aos/1056562464}.
\bibitem[{Garel(2001)}]{garel01}
\bibinfo{author}{Garel, B.} (\bibinfo{year}{2001}).
\newblock \bibinfo{title}{Likelihood ratio test for univariate gaussian
  mixtures}.
\newblock {\it \bibinfo{journal}{Journal of Statistical Planning and
  Inference}\/},  {\it \bibinfo{volume}{96}\/}, \bibinfo{pages}{325--350}.
\bibitem[{Garel(2005)}]{garel05}
\bibinfo{author}{Garel, B.} (\bibinfo{year}{2005}).
\newblock \bibinfo{title}{Asymptotic theory of the likelihood ratio test for
  the identification of a mixture}.
\newblock {\it \bibinfo{journal}{J. Statist. Plann. Inference}\/},  {\it
  \bibinfo{volume}{131}\/}, \bibinfo{pages}{271--296}. \URLprefix
  \url{https://doi.org/10.1016/j.jspi.2004.01.006}.
  \DOIprefix\doi{10.1016/j.jspi.2004.01.006}.
\bibitem[{Gassiat(2002)}]{gassiat02}
\bibinfo{author}{Gassiat, E.} (\bibinfo{year}{2002}).
\newblock \bibinfo{title}{Likelihood ratio inequalities with applications to
  various mixtures}.
\newblock (pp. \bibinfo{pages}{897--906}).
\newblock volume~\bibinfo{volume}{38}.
\newblock \URLprefix \url{https://doi.org/10.1016/S0246-0203(02)01125-1}.
  \DOIprefix\doi{10.1016/S0246-0203(02)01125-1} \bibinfo{note}{en l'honneur de
  J. Bretagnolle, D. Dacunha-Castelle, I. Ibragimov}.
\bibitem[{van~de Geer(2000)}]{geer00}
\bibinfo{author}{van~de Geer, S.~A.} (\bibinfo{year}{2000}).
\newblock {\it \bibinfo{title}{Empirical Processes in M-estimation}\/}
  volume~\bibinfo{volume}{6}.
\newblock \bibinfo{publisher}{Cambridge university press}.
\bibitem[{Ghosh \& Sen(1985)}]{ghoshsen85}
\bibinfo{author}{Ghosh, J.~K.}, \& \bibinfo{author}{Sen, P.~K.}
  (\bibinfo{year}{1985}).
\newblock \bibinfo{title}{On the asymptotic performance of the log likelihood
  ratio statistic for the mixture model and related results}.
\newblock In {\it \bibinfo{booktitle}{Proceedings of the {B}erkeley conference
  in honor of {J}erzy {N}eyman and {J}ack {K}iefer, {V}ol. {II} ({B}erkeley,
  {C}alif., 1983)}\/} Wadsworth Statist./Probab. Ser. (pp.
  \bibinfo{pages}{789--806}).
\newblock \bibinfo{publisher}{Wadsworth, Belmont, CA}.
\bibitem[{Hall \& Stewart(2005)}]{hallstewart05}
\bibinfo{author}{Hall, P.}, \& \bibinfo{author}{Stewart, M.}
  (\bibinfo{year}{2005}).
\newblock \bibinfo{title}{Theoretical analysis of power in a two-component
  normal mixture model}.
\newblock {\it \bibinfo{journal}{J. Statist. Plann. Inference}\/},  {\it
  \bibinfo{volume}{134}\/}, \bibinfo{pages}{158--179}.
\bibitem[{Hartigan(1985)}]{hartigan85}
\bibinfo{author}{Hartigan, J.~A.} (\bibinfo{year}{1985}).
\newblock \bibinfo{title}{A failure of likelihood asymptotics for normal
  mixtures}.
\newblock In {\it \bibinfo{booktitle}{Proceedings of the {B}erkeley conference
  in honor of {J}erzy {N}eyman and {J}ack {K}iefer, {V}ol. {II} ({B}erkeley,
  {C}alif., 1983)}\/} Wadsworth Statist./Probab. Ser. (pp.
  \bibinfo{pages}{807--810}).
\newblock \bibinfo{publisher}{Wadsworth, Belmont, CA}.
\bibitem[{Hui et~al.(2015)Hui, Warton \& Foster}]{hui15}
\bibinfo{author}{Hui, F. K.~C.}, \bibinfo{author}{Warton, D.~I.}, \&
  \bibinfo{author}{Foster, S.~D.} (\bibinfo{year}{2015}).
\newblock \bibinfo{title}{Order selection in finite mixture models: complete or
  observed likelihood information criteria?}
\newblock {\it \bibinfo{journal}{Biometrika}\/},  {\it
  \bibinfo{volume}{102}\/}, \bibinfo{pages}{724--730}. \URLprefix
  \url{https://doi.org/10.1093/biomet/asv027}.
  \DOIprefix\doi{10.1093/biomet/asv027}.
\bibitem[{H{\"u}sler(1990)}]{huesler90}
\bibinfo{author}{H{\"u}sler, J.} (\bibinfo{year}{1990}).
\newblock \bibinfo{title}{Extreme values and high boundary crossings of locally
  stationary {G}aussian processes}.
\newblock {\it \bibinfo{journal}{Ann. Probab.}\/},  {\it
  \bibinfo{volume}{18}\/}, \bibinfo{pages}{1141--1158}.
\bibitem[{H{\"u}sler(1995)}]{huesler95}
\bibinfo{author}{H{\"u}sler, J.} (\bibinfo{year}{1995}).
\newblock \bibinfo{title}{A note on extreme values of locally stationary
  {G}aussian processes}.
\newblock {\it \bibinfo{journal}{J. Statist. Plann. Inference}\/},  {\it
  \bibinfo{volume}{45}\/}, \bibinfo{pages}{203--213}.
\newblock \bibinfo{note}{{E}xtreme value theory and applications (Villeneuve
  d'Ascq, 1992)}.
\bibitem[{Ingster(1997)}]{ingster97}
\bibinfo{author}{Ingster, Y.~I.} (\bibinfo{year}{1997}).
\newblock \bibinfo{title}{Some problems of hypothesis testing leading to
  infinitely divisible distributions}.
\newblock {\it \bibinfo{journal}{Math. Methods Statist.}\/},  {\it
  \bibinfo{volume}{6}\/}, \bibinfo{pages}{47--69}.
\bibitem[{Ingster(2001)}]{ingster01}
\bibinfo{author}{Ingster, Y.~I.} (\bibinfo{year}{2001}).
\newblock \bibinfo{title}{Adaptive detection of a signal of growing dimension.
  {I}}.
\newblock {\it \bibinfo{journal}{Math. Methods Statist.}\/},  {\it
  \bibinfo{volume}{10}\/}, \bibinfo{pages}{395--421 (2002)}.
\newblock \bibinfo{note}{Meeting on Mathematical Statistics (Marseille, 2000)}.
\bibitem[{Ingster(2002)}]{ingster02}
\bibinfo{author}{Ingster, Y.~I.} (\bibinfo{year}{2002}).
\newblock \bibinfo{title}{Adaptive detection of a signal of growing dimension.
  {II}}.
\newblock {\it \bibinfo{journal}{Math. Methods Statist.}\/},  {\it
  \bibinfo{volume}{11}\/}, \bibinfo{pages}{37--68}.
\bibitem[{Leadbetter et~al.(1983)Leadbetter, Lindgren \&
  Rootz{\'e}n}]{leadbetter83}
\bibinfo{author}{Leadbetter, M.}, \bibinfo{author}{Lindgren, G.}, \&
  \bibinfo{author}{Rootz{\'e}n, H.} (\bibinfo{year}{1983}).
\newblock {\it \bibinfo{title}{Extremes and Related Properties of Random
  Sequences and Processes}\/}.
\newblock \bibinfo{publisher}{Springer-Verlag}.
\bibitem[{Lehmann \& Romano(2005)}]{lehmannromano05}
\bibinfo{author}{Lehmann, E.~L.}, \& \bibinfo{author}{Romano, J.~P.}
  (\bibinfo{year}{2005}).
\newblock {\it \bibinfo{title}{Testing statistical hypotheses}\/}.
\newblock Springer Texts in Statistics (\bibinfo{edition}{3rd} ed.).
\newblock \bibinfo{publisher}{Springer, New York}.
\bibitem[{Liu et~al.(2003)Liu, Pasarica \& Shao}]{liu03}
\bibinfo{author}{Liu, X.}, \bibinfo{author}{Pasarica, C.}, \&
  \bibinfo{author}{Shao, Y.} (\bibinfo{year}{2003}).
\newblock \bibinfo{title}{Testing homogeneity in gamma mixture models}.
\newblock {\it \bibinfo{journal}{Scand. J. Statist.}\/},  {\it
  \bibinfo{volume}{30}\/}, \bibinfo{pages}{227--239}.
\bibitem[{Liu \& Shao(2003)}]{liushao03}
\bibinfo{author}{Liu, X.}, \& \bibinfo{author}{Shao, Y.}
  (\bibinfo{year}{2003}).
\newblock \bibinfo{title}{Asymptotics for likelihood ratio tests under loss of
  identifiability}.
\newblock {\it \bibinfo{journal}{Annals of Statistics}\/},  {\it
  \bibinfo{volume}{31}\/}.
\bibitem[{Liu \& Shao(2004)}]{liushao04}
\bibinfo{author}{Liu, X.}, \& \bibinfo{author}{Shao, Y.}
  (\bibinfo{year}{2004}).
\newblock \bibinfo{title}{Asymptotics for the likelihood ratio test in a
  two-component normal mixture model}.
\newblock {\it \bibinfo{journal}{J. Statist. Plann. Inference}\/},  {\it
  \bibinfo{volume}{123}\/}, \bibinfo{pages}{61--81}.
\bibitem[{Piterbarg(1996)}]{piterbarg96}
\bibinfo{author}{Piterbarg, V.~I.} (\bibinfo{year}{1996}).
\newblock {\it \bibinfo{title}{Asymptotic Methods in the Theory of Gaussian
  Processes and Fields}\/}.
\newblock \bibinfo{publisher}{American Mathematical Society}.
\bibitem[{Porter \& Stewart(2020)}]{porterstewart20}
\bibinfo{author}{Porter, T.}, \& \bibinfo{author}{Stewart, M.}
  (\bibinfo{year}{2020}).
\newblock \bibinfo{title}{Beyond hc: More sensitive tests for rare/weak
  alternatives}.
\newblock {\it \bibinfo{journal}{Ann. Statist.}\/},  {\it
  \bibinfo{volume}{48}\/}, \bibinfo{pages}{2230--2252}.
  \DOIprefix\doi{10.1214/19-AOS1885}.
\bibitem[{Stewart \& Robinson(2003)}]{stewartrobinson03}
\bibinfo{author}{Stewart, M.}, \& \bibinfo{author}{Robinson, J.}
  (\bibinfo{year}{2003}).
\newblock \bibinfo{title}{Extremes of normed empirical moment generating
  function processes}.
\newblock {\it \bibinfo{journal}{Extremes}\/},  {\it \bibinfo{volume}{6}\/},
  \bibinfo{pages}{319--333 (2005)}.
\bibitem[{van~der Vaart \& Wellner(1996)}]{vanderwellner96}
\bibinfo{author}{van~der Vaart, A.~W.}, \& \bibinfo{author}{Wellner, J.~A.}
  (\bibinfo{year}{1996}).
\newblock {\it \bibinfo{title}{Weak Convergence and Empirical Processes}\/}.
\newblock Springer Series in Statistics.
\newblock \bibinfo{publisher}{Springer}.
\bibitem[{Yang(2005)}]{yang05}
\bibinfo{author}{Yang, Y.} (\bibinfo{year}{2005}).
\newblock \bibinfo{title}{Can the strengths of {AIC} and {BIC} be shared? {A}
  conflict between model indentification and regression estimation}.
\newblock {\it \bibinfo{journal}{Biometrika}\/},  {\it \bibinfo{volume}{92}\/},
  \bibinfo{pages}{937--950}. \URLprefix
  \url{http://biomet.oxfordjournals.org/content/92/4/937.abstract}.
  \DOIprefix\doi{10.1093/biomet/92.4.937}.
  \href{http://arxiv.org/abs/http://biomet.oxfordjournals.org/content/92/4/937.full.pdf+html}{\tt
  arXiv:http://biomet.oxfordjournals.org/content/92/4/937.full.pdf+html}.

\end{thebibliography}

\appendix
\section{R Code for  parametric bootstrap}
\label{sec:org889b3aa}
\end{document}